\documentclass{siamart1116}

\usepackage{amsopn}
\usepackage{amsmath}
\usepackage{amssymb}
\usepackage{afterpage}
\usepackage{bm}
\usepackage{booktabs}
\usepackage{multirow}
\usepackage[shortlabels]{enumitem}
\usepackage{algorithm,algorithmic}
\usepackage{datetime}
\usepackage[margin=1.38in,bottom=1.2in]{geometry}
\usepackage[protrusion=true,tracking=false,kerning=true]{microtype}
\usepackage{color}
\usepackage[framed,numbered,autolinebreaks]{mcode} 

\newcommand{\alfonso}{\texttt{alfonso}}
\renewcommand{\epsilon}{\varepsilon}
\newcommand{\defeq}{\ensuremath{\overset{\mathrm{def}}{=}}}
\newcommand{\deletethis}[1]{{}}

\newcommand{\vA}{\mathbf{A}}
\newcommand{\vb}{\mathbf{b}}
\newcommand{\vB}{\mathbf{B}}
\newcommand{\vc}{\mathbf{c}}
\newcommand{\vd}{\mathbf{d}}
\newcommand{\vdelta}{\bm{\delta}}
\newcommand{\vDelta}{\bm{\Delta}}
\newcommand{\ve}{\mathbf{e}}
\newcommand{\cE}{\mathcal{E}}
\newcommand{\cF}{\mathcal{F}}

\newcommand{\vG}{\mathbf{G}}
\newcommand{\vh}{\mathbf{h}}
\newcommand{\vH}{\mathbf{H}}
\newcommand{\vI}{\mathbf{I}}
\newcommand{\vM}{\mathbf{M}}
\newcommand{\vlambda}{{\bm{\lambda}}}
\newcommand{\vL}{\mathbf{L}}
\newcommand{\vmu}{{\bm{\mu}}}
\newcommand{\cN}{\mathcal{N}}
\newcommand{\Oh}{\mathcal{O}}
\newcommand{\vOmega}{\bm{\Omega}}
\newcommand{\vone}{\mathbf{1}}
\newcommand{\cP}{\mathcal{P}}
\newcommand{\vP}{\mathbf{P}}
\newcommand{\cQ}{\mathcal{Q}}
\newcommand{\R}{\mathbb{R}}
\newcommand{\Vs}{\mathbf{s}} 
\newcommand{\cS}{\mathcal{S}}
\newcommand{\vSigma}{\bm{\Sigma}}
\newcommand{\T}{{\!\mathrm{T}}}
\newcommand{\vt}{\mathbf{t}}
\newcommand{\vu}{\mathbf{u}}
\newcommand{\cV}{\mathcal{V}}

\newcommand{\vV}{\mathbf{V}}
\newcommand{\vx}{\mathbf{x}}
\newcommand{\vX}{\mathbf{X}}
\newcommand{\vy}{\mathbf{y}}
\newcommand{\vz}{\mathbf{z}}

\newcommand{\vzero}{\mathbf{0}}

\newcommand{\af}{\alpha}

\newcommand{\dt}{\Delta_{\tau}}
\newcommand{\dk}{\Delta_{\kappa}}
\newcommand{\dx}{\Delta_{\vx}}
\newcommand{\dy}{\Delta_{\vy}}
\newcommand{\ds}{\Delta_{\Vs}}
\newcommand{\dz}{\Delta_{\vz}}

\theoremstyle{plain}
\theoremstyle{definition}

\newtheorem{example}{Example}

\newcommand{\TheTitle}{alfonso: Matlab package for nonsymmetric conic optimization}
\newcommand{\TheAuthors}{D{\'a}vid PAPP and Sercan YILDIZ}

\headers{\TheTitle}{\TheAuthors}

\title{%
	{\TheTitle}%
	\thanks{%
		\textbf{Submitted to \emph{INFORMS Journal on Computing} on August 19, 2020. Find the published version on the journal's website.}
		\funding{DP: This material is based upon work supported by the National Science Foundation under Grant No. DMS-1719828 and Grant No. DMS-1847865. SY: The material was based upon work partially supported by the National Science Foundation under Grant No. DMS-1638521 to the Statistical and Applied Mathematical Sciences Institute.}
	}%
}

\author{
	D{\'a}vid PAPP\thanks{North Carolina State University, Department of Mathematics. Email: \email{dpapp@ncsu.edu}.}
	\and
	Sercan YILDIZ\thanks{Qontigo, Inc. Email: \email{syildiz@qontigo.com}.}
}

\ifpdf
\hypersetup{
	pdftitle={\TheTitle},
	pdfauthor={\TheAuthors}
}
\fi

\usepackage{natbib}
\bibpunct[, ]{(}{)}{,}{a}{}{,}%
\begin{document}

\maketitle

\begin{abstract}
	We present \alfonso, an open-source Matlab package for solving conic optimization problems over nonsymmetric convex cones. The implementation is based on the authors' corrected analysis of a primal-dual interior-point method of Skajaa and Ye. This method enables optimization over any convex cone as long as a logarithmically homogeneous self-concordant barrier is available for the cone or its dual. This includes many nonsymmetric cones, for example, hyperbolicity cones and their duals (such as  sum-of-squares cones), semidefinite and second-order cone representable cones, power cones, and the exponential cone.

Besides enabling the solution of problems which cannot be cast as optimization problems over a symmetric cone, algorithms for nonsymmetric conic optimization can also offer performance advantages for problems whose symmetric cone programming representation requires a large number of auxiliary variables or has a special structure that can be exploited in the barrier computation.


The worst-case iteration complexity of \alfonso{} is the best known for non-symmetric cone optimization: $\Oh(\sqrt{\nu}\log(1/\epsilon))$ iterations to reach an $\epsilon$-optimal solution, where $\nu$ is the barrier parameter of the barrier function used in the optimization.

\alfonso{} can be interfaced with a Matlab function (supplied by the user) that computes the gradient and the Hessian of a barrier function for the cone. For convenience, a simplified interface is also available to optimize over the direct product of cones for which a barrier function has already been built into the software. This interface can be easily extended to include new cones.

Both interfaces are illustrated using the toy example of solving linear programs in standard form. Additionally, the oracle interface and the efficiency of \alfonso{} are demonstrated using an optimal design of experiments problem in which the tailored barrier computation greatly decreases the solution time compared to using state-of-the-art off-the-shelf conic optimization software.
\end{abstract}



\section{Introduction}

We present \alfonso{}, an open-source, Octave-compatible, Matlab package for solving optimization problems over (not necessarily symmetric) convex cones. More precisely, \alfonso{} can be used to solve primal-dual pairs of optimization problems of the form\\

\noindent \begin{minipage}{.45\linewidth}
\begin{equation}\tag{P}\label{eq:P}
\begin{aligned}
&\underset{\vx\in\R^n}{\text{minimize}}\quad  &&\vc^\T \vx\\
                &\text{subject to}  &&\vA\vx=\vb\\
                &                   &&\vx\in K
\end{aligned}
\end{equation}
\end{minipage}%
\begin{minipage}{.45\linewidth}
\begin{equation}\tag{D}\label{eq:D}
\begin{aligned}
&\underset{\Vs\in\R^n,\,\vy\in\R^m}{\text{maximize}}\quad   &&\vb^\T \vy\\
                &\text{subject to}  &&\vA^\T \vy+\Vs=\vc\\
                &                   &&\Vs\in K^*,
\end{aligned}
\end{equation}
\end{minipage}

\vspace{1.5ex}

\noindent where $K$ is a full-dimensional, pointed, closed, convex cone and $K^*$ is its dual. The only additional assumption $K$ needs to satisfy is that \emph{an efficient algorithm to compute the gradient and the Hessian of some logarithmically homogeneous self-concordant barrier function} of $K$ is available. (See Sec.~\ref{sec:background}.) As an automatic generalization, it is also sufficient to have such a barrier function for only the dual cone $K^*$, since we can apply \alfonso{} to the dual problem. A feasible initial point is not required, only an initial point in the interior of $K$.

We may assume without loss of generality that $\operatorname{rank}(\vA) = m$. 
If $\operatorname{rank}(\vA) < m$, then depending on whether $\vb\in\operatorname{range}(\vA)$ or not, either the equality constraints in the primal problem are inconsistent or some of the equalities are redundant and can be removed.

The set of problems \alfonso{} can solve includes optimization over many nonsymmetric cones of great interest, for example, hyperbolicity cones of efficiently computable hyperbolic polynomials \citep{Renegar2004} and their duals, sum-of-squares cones \cite[Chapter 3]{BlekhermanParriloThomas2013}, $\mathcal{L}_p$ cones \citep{GlineurTerlaky2004} and other flavors of (generalized) power cones \citep{RoyXiao2018}, and the exponential cone \citep{Chares2009}.

To maximally take advantage of this level of generality, \alfonso{} can be interfaced directly with a function handle to a \emph{membership and barrier function oracle}, a Matlab function that computes whether a given point is in the interior of the cone, and for interior points computes the gradient and Hessian of an appropriate barrier function. For convenience, we have also created a simplified interface that allows the user to specify $K$ as the direct product of known cones for which \alfonso{} already has oracles implemented. Through this interface, \alfonso{} is easily extensible: the barrier function of any cone may be implemented and then added to the list of cones accepted by this interface, by adding only a few additional lines to \alfonso's source code.

To our knowledge there are very few alternative conic optimization software that offer this level of generality and extensibility. SCS \citep{scs} is a first-order, operator splitting method that can solve conic optimization problems over any cone that is easy to orthogonally project to, and currently supports the exponential and power cones in addition to symmetric cones. ECOS \citep{ecos} is a second-order cone programming software whose latest version also handles exponential cone constraints.
DDS \citep{KarimiTuncel2019} is a recent solver that aims at the same level of generality as \alfonso{} with an entirely different approach to domain definition. Hypatia is a recently announced solver written in Julia, based on a similar algorithm as \alfonso{}  \citep{CoeyKapelevichVielma2020}, but at the time of writing this paper, the code does not appear to be publicly available.\footnote{As of October 2020, Hypatia has been released.} In the commercial domain, Mosek 9 is capable of solving conic optimization problems with any combination of symmetric, exponential and power cone constraints \citep{mosek9}, but it is neither open-source nor extensible with new cones.

Besides enabling the solution of problems which cannot be cast as optimization over a symmetric cone, algorithms for nonsymmetric conic optimization can also offer performance advantages for problems that can be written as optimization problems over symmetric cones. This is the case, for example, when the equivalent representation as an optimization problem over a symmetric cone requires an extended formulation with a large number of auxiliary variables, or when the representation has some special structure that all-purpose optimization software often do not take advantage of, such as Hankel, Toeplitz, or low-rank structures in semidefinite programming.
An example of a family of optimization problems that greatly benefit from a non-symmetric cone optimization approach is \emph{sum-of-squares optimization}; this application was the initial motivation for the development of \alfonso{} \citep{PappYildiz2019}. In Section~\ref{sec:design} we demonstrate another application in which \alfonso{} is  several orders of magnitude more efficient than the straightforward semidefinite programming approach.
\deletethis{Sum-of-squares (SOS) cones are cones of polynomials instrumental in the  global optimization of polynomials \citep{Lasserre2015, BlekhermanParriloThomas2013}. Even though SOS cones are semidefinite representable \citep{Nesterov2000}, which means that sum-of-squares optimization problems can be cast as semidefinite programs and solved using off-the-shelf semidefinite programming solvers \citep{gloptipoly, sostools}, it was shown recently by the authors that optimizing directly over the (non-symmetric) SOS cone can be several orders of magnitude faster than the conventional semidefinite programming approach \citep{PappYildiz2019}. This application was the initial motivation for the development of \alfonso. In Section~\ref{sec:design} we demonstrate another application in which \alfonso{} is several orders of magnitude more efficient than the straightforward semidefinite programming approach.} Additional complex examples with extensive computational results comparing an earlier version of the code with state-of-the-art off-the-shelf interior-point solvers can be found in the authors' recent work on polynomial optimization \citep{PappYildiz2019, Papp2019}.

The implementation is based on an interior-point method (IPM) originally proposed by \cite{SkajaaYe2015} and subsequently improved by the authors \citep{PappYildiz2017corrigendum}. The iteration complexity of the method matches the iteration complexity of popular algorithms for symmetric cone optimization.

The source code can be found at \url{https://github.com/dpapp-github/alfonso}.

\subsection{Mathematical Background}\label{sec:background}

\begin{definition}[proper cone]
We say that a set $K\subseteq\R^n$ is a \emph{cone} provided that for every $\vx\in K$ and $\lambda \geq 0$ we also have $\lambda\vx \in K$. A cone is \emph{proper} if it satisfies all of the following: it is closed, convex, has a non-empty interior, and does not contain a line.
\end{definition}

Aside from closedness and convexity, the remaining assumptions on $K$ are essentially without loss of generality in the sense that every finite-dimensional closed convex optimization problem can be equivalently written in the form of \eqref{eq:P} with an appropriate choice of the full-row-rank matrix $\vA$, vectors $\vb, \vc$ and the proper cone $K$.

\begin{definition}[self-concordant barrier function, LHSCB]
Let $K^\circ$ denote the interior of $K$. A function $f\colon K^\circ\mapsto\R$ is a \emph{barrier} function if $f(\vx_i)\to\infty$ for every sequence $\vx_1,\vx_2,\dots$ of points $\vx_i \in K^\circ$ converging to a boundary point of $K$. A barrier is \emph{self-concordant} if it is convex, three times continuously differentiable, and if for every $\vx\in K^\circ$ and $\vh\in\R^n$ the inequality
$ |D^3 f(\vx)[\vh,\vh,\vh]| \leq 2D^2 f(\vx)[\vh,\vh]^{3/2} $
holds. The self-concordant barrier $f$ is called \emph{logarithmically homogeneous} if there exists a scalar $\nu$ such that for every $\vx\in K^\circ$ and $t>0$ we have
$ f(t\vx) = f(\vx) - \nu \ln t. $
The scalar $\nu$ is called the \emph{barrier parameter} of $K$. 
\end{definition}


As a shorthand, we say that $f$ is a \emph{$\nu$-LHSCB for $K$} if $f$ is a logarithmically homogeneous self-concordant barrier whose domain is the interior of a proper convex cone $K$ and if the barrier parameter of $f$ is $\nu$. For the interested reader, the monograph \citep{NesterovNemirovskii1994} provides a comprehensive treatment of LHSCBs. Renegar's treatment \citep{Renegar2001} of the subject is also excellent. Many fundamental cones in the application of convex optimization have known and easily computable LHSCBs.
\deletethis{
It is known that every proper convex cone $K$ has an LHSCB with $K^\circ$ as its domain \cite[Section 2.5]{NesterovNemirovskii1994}, but such barriers are not always available in closed form; moreover, even the membership problem can be NP-hard for a proper convex cone. However, many fundamental cones in the application of convex optimization have known and easily computable LHSCBs. The barrier parameter is typically not greater than the dimension of these cones. Additionally, LHSCBs can be easily constructed for the direct (Cartesian) products and intersections of cones with known LHSCBs.
}

\begin{example}\label{ex:LHSCBs}\normalfont
The following examples are proper convex cones with the additional property that either the cone or its dual has a known LHSCB with easily computable derivatives. Only the first three cones are symmetric, the remaining ones are not.
\begin{enumerate}[{1}.1]
\item The function $f(x) = -\ln x$ is an LHSCB for $\R_+$, and more generally
$f(\vx) = -\sum_{i=1}^n \ln x_i$ is an $n$-LHSCB for $\R_+^n$.

\item The function
\[f(\vx) = -\ln\left(x_0^2 - \sum_{i=1}^n x_i^2\right)\] is a $2$-LHSCB for the \emph{second-order cone}
\[ \cQ_{n+1} \defeq \{ (x_0,\dots,x_n)\,|\, x_0 \geq \|(x_1,\dots,x_n)\| \}. \]
Note that its barrier parameter is independent of the dimension $n$.

\item \label{ex:sdp} The function
\[f(\vX) = -\ln\det\vX\]
is an $n$-LHSCB for the cone of $n\times n$ positive semidefinite real symmetric matrices.

\item The \emph{exponential cone} is the three-dimensional cone
\[ \cE \defeq \operatorname{cl}\left(\{\vx\in\R_+^2\times\R\,|\, x_1 > x_2 e^{x_3/x_2} \}\right). \]
The function
\[ f(\vx) = -\ln(x_1) - \ln(x_2) - \ln(x_2 \ln(x_1/x_2) - x_3) \]
is a $3$-LHSCB for this cone \cite[Chapter 2]{Chares2009}.

\item \label{ex:powercone} Suppose $\vlambda = (\lambda_1,\dots,\lambda_n) \in \R^n$ satisfies $\lambda_i>0$ for each $i$ and $\sum_{i=1}^n\lambda_i = 1$. Then the \emph{(generalized) power cone} with \emph{signature}
$\vlambda$ is the convex cone defined as
\begin{equation}\label{eq:power-cone-def}
    \cP_\vlambda \defeq \left\{ (\vx,z)\in \R_+^n \times \R \,\middle|\, |z|\leq \prod_{i=1}^n x_i^{\lambda_i} \right\}.
\end{equation}
The function
\[ f(\vx,z) = -\ln\left(\prod_{i=1}^n x_i^{2\lambda_i}-z^2\right) - \sum_{i=1}^n (1-\lambda_i)\ln(x_i) \] 
an $(n+1)$-LHSCB for this cone. This was first proven by \cite{RoyXiao2018}, who also study a number of related cones. 
The dual cone of $\cP_\vlambda$ is identical to the cone known in algebraic geometry as the \emph{SONC}, short for \emph{sum of nonnegative circuit polynomials} cone. \citep{IlimanDeWolff2016}.

\item A homogeneous $n$-variate polynomial $h$ of degree $d$ is said to be \emph{hyperbolic with respect to the point $\ve\in\R^n$} if $h(\ve) > 0$ and if for every $\vx\in\R^n$, the univariate polynomial $h(\vx+t\ve)$ has only real roots. 
The corresponding \emph{hyperbolicity cone} is the set
 \[\Lambda^{+}_{h,\ve} \defeq \{ \vx\in\R^n\,|\,h(\vx+t\ve)>0\;\;\forall\,t>0\}. \]
It can be shown that $\Lambda^{+}_{h,\ve}$ is a proper convex cone for which $-\ln h(\cdot)$ is a $d$-LHSCB \citep{Guler1997}. Because the determinant is a hyperbolic polynomial (with respect to the identity matrix) whose hyperbolicity cone is the semidefinite cone, optimization over hyperbolicity cones is a generalization of semidefinite programming.

\deletethis{
\item The \emph{spectral norm cone} $\cN_{m,n}$ is the epigraph of the spectral norm:
\[  \cN_{m,n} \defeq \left\{ (t,\vA)\in\R\times\R^{m \times n}\;\middle|\; t \geq \|\vA\|_1 \right\}. \]
Since $t \geq \|\vA\|_1$ if and only if $\left(\begin{smallmatrix}t\vI_m & \vA\\ \vA^\T & t\vI_n\end{smallmatrix}\right) \succcurlyeq \vzero$, Example \ref{item:Ax} below implies that
\[ f(t,\vA) = -\ln\det \left(\begin{smallmatrix}t\vI_m & \vA\\ \vA^\T & t\vI_n\end{smallmatrix}\right)\] 
is an $(m+n)$-LHSCB for this cone. It can be shown that the cone also admits a $\nu$-LHSCB with parameter $\nu=\min\{m,n\}+1$; see Proposition 5.4.6 in \citep{NesterovNemirovskii1994}.
}

\deletethis{
\item TOO MUCH TO ADD? Consider the set of $n\times n$ real symmetric matrices whose $2 \times 2$ (not necessarily leading) principal submatrices are all positive semidefinite
\begin{equation}\label{eq:SDSOS-dual}
\cS_n \defeq \left\{ \vX = (x_{ij})\in\R^{n \times n}\;\middle|\; 
\left(\begin{smallmatrix} x_{ii} & x_{ij}\\ x_{ji} & x_{jj} \end{smallmatrix}\right) \succcurlyeq \vzero \quad \forall 1\leq i<j\leq n \right\}.
\end{equation}
Its dual cone $\cS^*_n$ is known as the set of \emph{scaled diagonally dominant} (or \emph{SDSOS} for short) matrices, and has found applications in polynomial optimization \citep{AhmadiMajumdar2014, KuangGhaddarNaoumSawayaZuluaga2019} and power systems engineering \citep{KuangGhaddarNaoumSawayaZuluaga2017}. While we are not aware of any useful LHSCBs for the SDSOS cone, the representation of its dual \eqref{eq:SDSOS-dual} together with Example \ref{item:intersection} below imply that the $\cS_n$ admits an $\nu$-LHSCB with easily computable derivatives for $\nu=n(n-1)$. 
}

\deletethis{
\item \label{ex:SOS} Let $\cV_{n,d}$ denote the space $n$-variate of polynomials of degree $d$. We say that $p\in \cV_{n,2d}$ is \emph{sum-of-squares (SOS)} if there exist polynomials $q_1,\dots,q_M \in \cV_{n,d}$ satisfying $p = \sum_{i=1}^M q_i^2$. 
SOS cones are instrumental in polynomial optimization; see, e.g.  \citep{BlekhermanParriloThomas2013} and \citep{Laurent2009}.

The set $\Sigma_{n,2d}$ of SOS polynomials is a proper cone in $\cV_{n,2d}$ \citep[Thm.~17.1]{Nesterov2000}. While no practical LHSCB is known for $\Sigma_{n,2d}$, its dual $\Sigma_{n,2d}^*$ admits a $\nu$-LHSCB for $\nu = \dim(\cV_{n,d}) = \binom{n+d}{d}$ with easily computable derivatives, because there exists a linear map $\Lambda$ from $\operatorname{span}(\Sigma_{n,2d}^*)$ to the space of $\nu\times\nu$ real symmetric matrices such that
$(\Sigma_{n,2d}^*)^\circ = \{\vx\,|\,\Lambda(\vx) \succ 0\}$ \citep[Thm.~17.8]{Nesterov2000}. Therefore, combining Example~\ref{ex:sdp} and Example~\ref{item:Ax} below, the function
\[ f(\vx) = -\ln\det\Lambda(\vx) \]
is a $\nu$-LHSCB for $\Sigma_{n,2d}^*$. The computation of the barrier gradient and Hessian is particularly efficient if the polynomials are represented in an interpolant basis; see \citep[Theorem 3.5]{PappYildiz2019}.
}

\item \label{item:intersection} If $K_1,\dots,K_k$ are proper convex cones in $\R^n$ whose interiors have a non-empty intersection and $f_i$ is a $\nu_i$-LHSCB for $K_i$ ($i=1,\dots,k$), then $\sum_{i=1}^k f_i$ is a $\nu$-LHSCB for the intersection $\bigcap_{i=1}^k K_i$, with barrier parameter $\nu = \sum_{i=1}^k \nu_i$.

\item If $K_1\subseteq\R^{n_1}, \dots, K_k\subseteq\R^{n_k}$ are proper convex cones and $f_i$ is a $\nu_i$-LHSCB for $K_i$ ($i=1,\dots,k$), then $\sum_{i=1}^k f_i$ is a $\nu$-LHSCB for the product cone $K_1\times \cdots \times K_k$, with barrier parameter $\nu = \sum_{i=1}^k \nu_i$. 

\item \label{item:Ax} Let $K\subseteq\R^n$ be a proper convex cone with an LHSCB $f$, and $L$ be a linear subspace of $\R^n$ that intersects $K^\circ$. Then $f|_L$, the restriction of $f$ to the subspace $L$, is an LHSCB for $K\cap L$. Denoting the orthogonal projection matrix onto $L$ by $\vP_L$, the gradient and Hessian of this barrier function are $\vP_L \nabla f(\vx)$ and $\vP_L \nabla^2 f(\vx) \vP_L$, respectively.

In a different notation, if $K\subseteq\R^m$ is a proper convex cone with an LHSCB $f_K$ and $\vA\in \R^{m \times n}$ is a matrix whose range space intersects $K^\circ$, then the cone
\vspace{-.5ex}%
\[ C = \{ \vx \in \R^n \,|\, \vA\vx \in K \} \]
\vspace{-.5ex}%
is a proper convex cone,
and $f_C(\vx) = f_K(\vA\vx)$ is an LHSCB for $C$ whose gradient and Hessian are easily computable from the gradient and Hessian of $f_K$. A notable special case is when $K$ is the positive semidefinite cone; sets $C$ that can be written in this form are called \emph{spectrahedral cones} \citep[Chapter 2]{BlekhermanParriloThomas2013}.

Examples of spectrahedral cones which benefit from a non-symmetric cone optimization approach include the epigraph of the spectral norm \citep[Prop.~5.4.6]{NesterovNemirovskii1994}, also known as the \emph{spectral norm cone},
and the cone of \emph{sum-of-squares polynomials}, for which the LHSCB inherited from semidefinite programming is particularly efficiently computable when the polynomials are represented in an interpolant basis \citep{PappYildiz2019}.
\end{enumerate}
\end{example}

\smallskip

Example \ref{item:Ax} is particularly notable, because even though we have easily computable LHSCBs for these cones, there does not appear to be any straightforward way to construct easily computable barrier functions for their dual cones. Additional techniques to construct LHSCBs for convex cones from known LHSCBs of simpler cones can be found in \citep[Chapter 5]{NesterovNemirovskii1994}.

\deletethis{
There is one more fundamental property of LHSCBs we must review; it plays an important role in the initialization of the algorithm:
\begin{proposition}\label{thm:bijection}
For every $\vx\in K^\circ$ and every LHSCB $f$ of $K$, the vector $-\nabla f(\vx)$ belongs to $(K^*)^\circ$.
\end{proposition}


To compute an approximately optimal solution to \eqref{eq:P}-\eqref{eq:D} or a certificate of primal or dual infeasibility, \alfonso{} computes an $\varepsilon$-feasible solution to the so-called \emph{homogeneous self-dual embedding} of the problems \eqref{eq:P}-\eqref{eq:D}, derived from the optimality conditions of these problems with the introduction of two new (nonnegative) scalar variables, $\tau$ and $\kappa$. This homogeneous self-dual embedding is customarily written as follows:
\begin{gather*}
\begin{aligned}
                &           &\quad \vA  &\vx    & -\vb  & \tau  &&          &&          &&=\vzero\\
-\vA^\T         &\vy        &           &       & +\vc  & \tau  &&-\Vs      &&          &&=\vzero\\
\quad \vb^\T    &\vy        &-\vc^\T    &\vx    &       &       &&          &&-\kappa   &&=0
\end{aligned}\\
\begin{gathered}
\vy\in\R^k,\quad (\vx,\tau)\in K\times\R_+,\quad (\Vs,\kappa)\in K^*\times\R_+.\notag
\end{gathered}
\end{gather*}
In more concise notation, letting
\[\vG\defeq\left(
\begin{smallmatrix}
\vzero      & \vA       & -\vb\\
-\vA^\T     & \vzero    &  \vc\\
\vb^\T      &-\vc^\T    &  0
\end{smallmatrix}\right) \quad \text{ and } \quad \cF \defeq \R^k \times K \times \R_+ \times K^* \times \R_+, \]
the homogeneous self-dual embedding can be written as
\begin{equation}\label{eq:embed}
\vG\left(\begin{smallmatrix}\vy\\ \vx\\ \tau  \end{smallmatrix}\right) - \left(\begin{smallmatrix}\vzero\\ \Vs\\ \kappa \end{smallmatrix}\right) =
\left(\begin{smallmatrix}\vzero\\ \vzero\\ 0 \end{smallmatrix}\right)
\quad \text{ and } \quad (\vy,\vx,\tau,\Vs,\kappa) \in \cF.
\end{equation}
The following result is standard in interior-point methods; see, e.g., \cite[Lemma 1]{SkajaaYe2015}:
\begin{theorem}\label{thm:stopping}
Let $(\vy,\vx,\tau,\Vs,\kappa)$ be a feasible solution to the system \eqref{eq:embed}. Then $\vx^\T\Vs + \tau\kappa = 0$ and consequently $\tau\kappa = 0$.
\begin{itemize}
\item If $\kappa = 0 < \tau$, then $\vx^* \defeq \vx/\tau$ is an optimal solution to \eqref{eq:P} and $(\vy^*,\Vs^*) \defeq (\vy/\tau, \Vs/\tau)$ is an optimal solution \eqref{eq:D}.
\item If $\tau = 0 < \kappa$ and $\vb^\T\vy > 0$, then the primal problem \eqref{eq:P} is infeasible.
\item If $\tau = 0 < \kappa$ and $\vc^\T\vx < 0$, then the dual problem \eqref{eq:D} is infeasible.
\end{itemize}
\end{theorem}

The last two cases are not mutually exclusive.

In the rest of this section, in order to lighten the notation, we will frequently use the shorthand $\vz=(\vy,\vx,\tau,\Vs,\kappa)$ to denote the entire variable vector.
Motivated by Theorem \ref{thm:stopping}, we measure the progress of \alfonso{} by the violation of the linear constraints in the system \eqref{eq:embed} and by the \emph{complementarity gap}
\begin{equation}\label{eq:mu}
\mu(\vz) = \frac{\vx^\T\Vs+\tau\kappa}{\nu+1}.
\end{equation}

We compute an approximate solution to \eqref{eq:embed} by approximately following the central path associated with this system, given below. 
Given a scalar $\theta\geq 0$, the \emph{$\theta$-neighborhood of the central path} is the set of points
\[ \cN(\theta) \defeq \left\{ \vz=(\vy,\vx,\tau,\Vs,\kappa) \in \cF^\circ\;\middle|\; \left\|\left(\begin{smallmatrix} \nabla^2 f(\vx) & \vzero \\ \vzero & 1/\tau^2\end{smallmatrix}\right)^{-1/2}\psi(\vz,\mu(\vz))\right\| \leq \theta \mu(\vz) \right\},\]
where $\mu$ is defined in \eqref{eq:mu} and 
\begin{equation}\label{eq:psi}
\psi(\vz,t) = \left(\begin{smallmatrix}\Vs\\ \kappa\end{smallmatrix}\right) + t \left(\begin{smallmatrix}\nabla f(\vx)\\ -1/\tau \end{smallmatrix}\right).
\end{equation}
The algorithm belongs to the family of predictor-corrector methods. It alternates between \emph{predictor steps} that make progress towards satisfying the system \eqref{eq:embed} at the cost of increasing the radius $\theta$ of the neighborhood containing the current iterate and \emph{corrector steps} that bring the current iterate back to the vicinity of the central path while maintaining the previous progress. Both the predictor and the corrector steps are damped Newton iterations. In the predictor steps Newton's method is applied to \eqref{eq:embed} itself; in the corrector steps it is applied to \eqref{eq:embed} with a modified right-hand side. For efficiency, the code includes both hard-wired short step sizes that guarantee that the iterates remain in the required neighborhoods of the central path and a line search procedure that searches for the (approximately) largest acceptable step size in each predictor step.

Using the self-dual embedding and the above definition of the neighborhood of the central path, initialization is straightforward. For every $\vx^0 \in K^\circ$, the vector
\begin{equation}\label{eq:z0}
\vz^0 = (\vy^0,\vx^0,\tau^0,\Vs^0,\kappa^0) = (\vzero,\vx^0,1,-\nabla f(\vx^0), 1)
\end{equation} is on the central path, that is, $\vz^0 \in \cN(0)$. In particular, Proposition \ref{thm:bijection} shows that $\Vs^0\in (K^*)^\circ$, therefore $\vz^0\in\cF^\circ$. From \eqref{eq:gTx} we also see that this initial point has complementarity gap $\mu(\vz^0) = 1$.

The pseudocode for the complete algorithm is listed in Algorithm~\ref{alg:SY}.

\begin{figure*}
{
\centering
	\begin{minipage}[tbh]{\linewidth}
		\begin{algorithm}[H]
			\caption{\small \texttt{alfonso}: ALgorithm for Non-Symmetric conic Optimization}
			\label{alg:SY}
			\begin{algorithmic}
			\small
				\STATE \hspace{-1em}\textbf{Parameters:} Radius $0<\eta<1$, step sizes $\af_p$ and $\af_c$, and integer $r_c>0$ chosen following \citep{PappYildiz2017corrigendum}.
				\STATE \hspace{-1em}\textbf{Input:} Problem data $(\vA, \vb, \vc$), an LHSCB $f$ for $K$, and an initial solution $\vz^0 = (\vy^0, \vx^0, \tau^0, \Vs^0, \kappa^0) \in \cN(\eta)$.
				\LOOP
				\STATE \textbf{Termination?}
				\STATE If the termination criterion is satisfied, stop and return $\vz$.
				\STATE \textbf{Predictor step}
				\STATE Compute the matrix $\vH = \left(\begin{smallmatrix} \nabla^2 f(\vx) & \vzero \\ \vzero & 1/\tau^2\end{smallmatrix}\right)$, and solve for $\dz=(\dy,\dx,\dt,\ds,\dk)$ the system
				 \begin{equation*}
					\begin{aligned}
					 \vG\begin{pmatrix}\dy\\ \dx \\ \dt\end{pmatrix}-\begin{pmatrix}\vzero\\ \ds\\ \dk\end{pmatrix}&=-\vG\begin{pmatrix}\vy\\ \vx \\ \tau\end{pmatrix}+\begin{pmatrix}\vzero\\ \Vs\\ \kappa\end{pmatrix},\\
					 \begin{pmatrix}\ds\\ \dk\end{pmatrix}+\mu(\vz)\vH \begin{pmatrix}\dx\\ \dt\end{pmatrix}&=-\begin{pmatrix}\Vs\\ \kappa\end{pmatrix}.
					\end{aligned}
				\end{equation*}\vspace{-1em}
				\STATE Set $\vz\leftarrow \vz+\af_p\dz$.
				\STATE \textbf{Corrector steps}
				\FOR{$i=1,\ldots,r_c$}
				\STATE Compute the vector $\psi(\vz,\mu(\vz))$ from \eqref{eq:psi} and the matrix $\vH = \left(\begin{smallmatrix} \nabla^2 f(\vx) & \vzero \\ \vzero & 1/\tau^2\end{smallmatrix}\right)$, and solve for $\dz=(\dy,\dx,\dt,\ds,\dk)$ the system
				 \begin{equation*}
					\begin{aligned}
					 \vG\begin{pmatrix}\dy\\ \dx \\ \dt\end{pmatrix}-\begin{pmatrix}\vzero\\ \ds\\ \dk\end{pmatrix}&=\vzero,\\
					 \begin{pmatrix}\ds\\ \dk\end{pmatrix}+\mu(\vz)\vH \begin{pmatrix}\dx\\ \dt\end{pmatrix}&=-\psi(\vz,\mu(\vz)).
					\end{aligned}
				\end{equation*}\vspace{-1em}
				\STATE Set $\vz\leftarrow \vz+\af_c\dz$.
				\ENDFOR
				\ENDLOOP
			\end{algorithmic}
		\end{algorithm}
	\end{minipage}
	\vspace{18pt}
}
\end{figure*}

We now state the main convergence and complexity result from our analysis. 

\begin{proposition}\label{thm:alg-complex} \textup{\citep{PappYildiz2017corrigendum}}
For every $\varepsilon>0$, the parameters $\eta$, $\af_p$, $\af_c$, and $r_c$ can be chosen such that, given any initial solution $\vz^0 \in\cN(\eta)$, Algorithm~\ref{alg:SY} terminates with a solution $\vz^* = (\vy^*,\vx^*,\tau^*,\Vs^*,\kappa^*) \in\cN(\eta)$ that satisfies
\begin{equation}\label{eq:alg-complex}
\mu(\vz^*)\leq\varepsilon\mu(\vz^0)\quad\text{and}\quad
\left\|\vG\left(\begin{smallmatrix}\vy^*\\ \vx^*\\ \tau^*  \end{smallmatrix}\right) - \left(\begin{smallmatrix}\vzero\\ \Vs^*\\ \kappa^* \end{smallmatrix}\right)\right\|\leq\varepsilon\left\|\vG\left(\begin{smallmatrix}\vy^0\\ \vx^0\\ \tau^0  \end{smallmatrix}\right) - \left(\begin{smallmatrix}\vzero\\ \Vs^0\\ \kappa^0 \end{smallmatrix}\right)\right\|
\end{equation}
in $\Oh(\sqrt{\nu}\log(1/\epsilon))$ iterations.
\end{proposition}

Recall from above that every user-supplied $\vx^0\in K^\circ$ can be extended to an initial $\vz^0 \in \cN(0)$ with $\mu(\vz^0) = 1$. A solution satisfying \eqref{eq:alg-complex} is referred to as an \emph{$\epsilon$-approximate} solution. Having found an $\epsilon$-approximate solution, Theorem~\ref{thm:stopping} is used to (numerically) identify infeasible or unbounded instances, or return a near-optimal solution.

Note that in Algorithm~\ref{alg:SY} and Proposition~\ref{thm:alg-complex} the predictor step size $\af_p$ is treated as a fixed parameter. The analysis of the algorithm and the iteration complexity result can be extended in a straightforward manner to a variant that instead uses line search.
The details of the analysis can be found in \citep{PappYildiz2017corrigendum}.
}

\subsection{The Algorithm and its Complexity}
 
The implementation is based on an IPM applied to a homogeneous self-dual embedding of \eqref{eq:P}-\eqref{eq:D} that was originally proposed by \cite{SkajaaYe2015} and subsequently improved by the authors. We refer the reader to \citep{PappYildiz2017corrigendum} for the details of the algorithm and its analysis and \citep[Sec.~2]{PappYildiz2019} for a brief summary.

This method is one of the theoretically most efficient algorithms applicable to nonsymmetric cone programming; its worst-case iteration complexity matches the iteration complexity of successful IPMs for symmetric cones. The main convergence and complexity result from our analysis \citep[Prop.~2.1]{PappYildiz2019} can be summarized as follows: the number of iterations and number of calls to the membership and barrier function oracle required to reduce the primal and dual infeasibility and complementarity metrics to $\epsilon$ times their initial value are $\Oh(\sqrt{\nu}\log(1/\epsilon))$ where $\nu$ is the barrier parameter of the barrier function. For most cones this means $\Oh(\sqrt{n}\log(1/\epsilon))$ iterations and oracle calls.

\section{Interfaces}

\subsection{Installation} \alfonso{} is entirely written in Matlab m-code, and is thus portable and easy to install: unzip the downloaded files in any directory and add the \texttt{src} subdirectory of the package to the Matlab (or Octave) path. (One of the examples not detailed in this paper requires additional packages.)

\subsection{Input interfaces}

An instance of the optimization problem \eqref{eq:P} can be described by the problem data $(\vA,\vb,\vc)$ and the cone $K$. Because of the level of generality \alfonso{} is aimed at, there are two ways to specify the cone when interfacing with the code.

\subsubsection{The oracle interface}
The cone $K$ can be specified using a \emph{membership and barrier function oracle}, which is a subroutine with the following signature:
\begin{lstlisting}[aboveskip=2ex,belowskip=1ex]
function [in, g, H, L] = oracle(x, bParams)
\end{lstlisting}
The first input argument \mcode{x} represents the primal vector $\vx$, which is the oracle's input.

The second argument \mcode{bParams} is an optional one that can be used to specify other parameters for the barrier function. For example, if \mcode{oracle} implements an LHSCB for the generalized power cone with signature $\vlambda$ (recall Example \ref{ex:powercone}), then it is convenient to pass $\vlambda$ as a parameter. If necessary, multiple parameters that cannot be conveniently passed as a single vector can be passed using a \mcode{struct} for \mcode{bParams}.

If $f$ denotes the LHSCB implemented in \mcode{oracle}, then the four outputs of the oracle are:
\begin{itemize}
	\item \mcode{in}: a Boolean flag that is \mcode{true} if $\vx \in K^\circ$ and \mcode{false} otherwise.
	\item \mcode{g}: a vector whose value is the gradient $\nabla f(\vx)$ if $\vx\in K^\circ$. Its value is ignored otherwise.
	\item \mcode{H}: a matrix whose value is the Hessian $\nabla^2 f(\vx)$ if $\vx\in K^\circ$. Its value is ignored otherwise.
	\item \mcode{L}: a lower triangular Cholesky factor of the Hessian.
\end{itemize}

\alfonso{} frequently calls the oracle with only the first or the first two output arguments. Unless all output parameters can be computed very efficiently, it is highly recommended that the oracle only computes the necessary output arguments, using Matlab's \mcode{nargout} feature.

If $\vH$ or $\vL$ is sparse, they should be computed as \mcode{sparse} matrices. The Cholesky factor can often be determined in closed form; otherwise one may always resort to the following generic code snippet to compute $\vL$ from $\vH$:
\begin{lstlisting}[aboveskip=2ex,belowskip=1ex]
if nargout > 3
	[L,err] = chol(H,'lower');
	if err > 0
		in = false; g = NaN; H = NaN; L = NaN;
		return;
	end
end
\end{lstlisting}

Lastly, \alfonso{} needs a starting point for the optimization. Only a primal initial point is needed, in the interior of $K$; \alfonso{} automatically computes an initial primal-dual iterate on the central path.
\deletethis{As discussed above, any point in the interior of $K$ suffices; more central points generally work better. \alfonso{} always uses $\vz^0$ defined in \eqref{eq:z0} for the initial primal-dual iterate, computed from the user-supplied $\vx^0$.}

Having all of this ready, the optimization problem can be solved by calling 
\begin{lstlisting}[aboveskip=2ex,belowskip=1ex]
alfonso(probData, x0, @oracle, bParams, opts)
\end{lstlisting}
The first argument is a Matlab \mcode{struct} with three mandatory fields, \mcode{A}, \mcode{b}, and \mcode{c}, and it contains the problem data. The second argument is the initial point. The third is a function handle to the membership and barrier function oracle, while the fourth (optional) argument is the parameter to be passed to the oracle as its second argument.

The last optional argument \mcode{opts} is a structure specifying the optimization options. See Section~\ref{sec:options} for more details on algorithmic and other options, and Section~\ref{sec:LP_example_oracle} for a complete example of how an optimization problem can be set up and solved using this interface.

\subsubsection{The simple interface}
The goal of the simple interface is to facilitate the reuse of previously implemented barrier functions. In the simple interface the cone is specified as a Cartesian product $K_1\times \cdots \times K_k$ of known cones $K_i$, passed to \alfonso{} as a Matlab cell array of structures whose $i$th element describes $K_i$.

In the cone array \mcode{K}, each element \mcode{K\{i\}} has two mandatory fields: \mcode{K\{i\}.type}, a string that specifies the cone $K_i$, and \mcode{K\{i\}.dim}, a string that specifies the dimension of the cone. 
The already built-in cones include:
\begin{itemize}
	\item \mcode{type = 'l'} or \mcode{'lp'} represents the nonnegative orthant.
	\item \mcode{type = 'soc'} or \mcode{'socp'} represents the second-order cone.
	\item \mcode{type = 'exp'} represents the exponential cone.
	\item \mcode{type = 'gpow'} represents a generalized power cone (defined in Example~\ref{ex:powercone}). The parameter $\vlambda$ must be specified in the field $\mcode{K\{i\}.lambda}$ as an additional vector.
\end{itemize}

Deviating slightly from the theory, variables in \alfonso{} are allowed to be \emph{free}, that is, not to be a member of any cone. This can be specified using \mcode{K\{i\}.type = 'free'}. Free variables are handled by placing them in a second-order cone using a single additional dummy variable, which is a common strategy in conic optimization attributed to \cite{Andersen2002}, and is also used, for instance, in SeDuMi. 

For example, the cell array
\begin{lstlisting}[aboveskip=2ex,belowskip=1ex]
K{1}.type = 'socp';    % second-order cone
K{1}.dim  = 10;
K{2}.type = 'free';    % free variables
K{2}.dim  = 6;
K{3}.type = 'lp';      % nonnegative orthant
K{3}.dim  = 10;
K{4}.type = 'exp';     % exponential cone, always 3-dimensional
\end{lstlisting}
defines the cone 
$K = \cQ_{10} \times \R^6 \times \R_+^{10} \times \cE$.

When $K$ is the Cartesian product of known cones, it is not necessary to provide an initial point; \alfonso{} defaults to the concatenation of known, ``central'', interior points of these cones. The syntax of the simple interface is
\begin{lstlisting}[aboveskip=2ex,belowskip=1ex]
alfonso_simple(c, A, b, K, x0, opts)
\end{lstlisting}
where the first four arguments are as described above, \mcode{x0} is the optional initial point (that can be set to \mcode{[]} for the default value) and the also optional \mcode{opts} argument is the same options structure as used in the oracle interface. (See Section~\ref{sec:options} for more details on the options.)

\subsection{Outputs}

Regardless of which interface is used, \alfonso{} returns a single structure as a result with over 20 fields that contain various diagnostic elements and information about the optimization process in addition to the primal and dual solutions. The comments in the header of \texttt{alfonso.m} contain a detailed description of all of them; here we only summarize the most important ones:
\begin{itemize}
\item \mcode{status}: an integer representing the solver status when the solver stopped. Its value is \mcode{1} if an approximately optimal solution was found.
\item \mcode{statusString}: the same information as \mcode{status} but in a human-readable format.
\item \mcode{x}, \mcode{s}, and \mcode{y}: the final primal and dual iterates.
\item \mcode{pObj} and \mcode{dObj}: final primal and dual objective function values.
\item \mcode{time}: the solution (wall-clock) time in seconds.
\end{itemize}

\subsection{Algorithmic and Other Options}\label{sec:options}


Options for \alfonso{} can be set using the optional last argument to the \mcode{alfonso()} or \mcode{alfonso_simple()} function. This argument is a structure (\mcode{struct}) with fields set to their desired values. Any options not specified this way will take their default values, which are detailed in the header of \texttt{alfonso.m}. The options that the users are most likely to want to change are the following:
\begin{itemize}
	\item \mcode{optimTol}: optimality tolerance $\epsilon$. 
	Default value: \mcode{1e-6}.
	\item \mcode{verbose}: a Boolean flag controlling the output level. \deletethis{If it is set to \mcode{0} (\mcode{false}), then \alfonso{} produces no output; if it is set to \mcode{1} (\mcode{true}), then the progress of the algorithm is printed in each iteration.}Default value: \mcode{1}. 	
\end{itemize}

The remaining options adjust various parameters of the algorithm (such as the line search procedure); these are documented in the header of \texttt{alfonso.m} and are omitted here, as changing them is only recommended in very particular situations.

\subsection{A minimal example: solving linear programs}\label{sec:LP_example_oracle}


In this section, we use the toy example of solving linear programs in standard form to illustrate how problem data is structured for each of the two interfaces. This example (with additional comments) is also included in the package in the files \texttt{random\_lp.m} and \texttt{random\_lp\_simple.m} in the directory \texttt{examples/random\_lp}. Additional examples can be found in the \texttt{examples} subdirectory of the code.


\subsubsection{The oracle interface}
To solve a linear program using the oracle interface, the user must implement a Matlab function that solves the membership problem and (for points in the interior) computes the gradient and factors the Hessian of an LHSCB for the nonnegative orthant $K=\R_+^n$. For the nonnegative orthant we use the logarithmic barrier $f$ given by $f(\vx) = -\sum_{i=1}^n \ln(x_i)$. A straightforward implementation is shown on Fig.~\ref{fig:gH_LP}. For efficiency, we use sparse matrices. The second input argument of the barrier function (that allows the passing of parameters) is not used.

\begin{figure}[!b]
\caption{A membership and barrier function oracle for solving linear programs in standard form.}\label{fig:gH_LP}
\begin{lstlisting}[aboveskip=1ex,belowskip=0ex]
function [in, g, H, L] = gH_lp(x, ~)
    n  = length(x);
    in = min(x)>0;
    if in
        g = -1./x;
        H = sparse(1:n,1:n,x.^(-2),n,n,n);
        L = sparse(1:n,1:n,-g,n,n,n);
    else
        g = NaN; H = NaN; L = NaN;
    end
end
\end{lstlisting}
\end{figure}

With the oracle \mcode{gH_lp()} ready, a linear program in standard form, with problem data $\vA$, $\vb$, and $\vc$ as in \eqref{eq:P}, can be solved by simply calling
\begin{lstlisting}[aboveskip=2ex,belowskip=1ex]
probData = struct('c', c, 'A', A, 'b', b);
results = alfonso(probData, x0, @gH_lp);
\end{lstlisting}
where \mcode{x0} is any componentwise positive initial point, e.g., the all-ones vector \mcode{ones(n,1)}. The optimal solution will be returned in \mcode{results.x}.

If any options are to be changed, the second line needs to include the options structure. In the following example, we decrease the optimality tolerance:
\begin{lstlisting}[aboveskip=2ex,belowskip=1ex]
opts.optimTol = 1e-7;
results = alfonso(probData, x0, @gH_lp, [], opts);
\end{lstlisting}
The empty list in the fourth argument is a placeholder for the optional parameters to pass to the function \mcode{gH_lp}, which is not used in this example.

\subsubsection{The simple interface} Using the simple interface, the user only needs to represent the $n$-dimensional nonnegative orthant in a cone structure (cell array) as follows:
\begin{lstlisting}[aboveskip=2ex,belowskip=1ex]
K{1} = struct('type', 'lp', 'dim', n);
results = alfonso_simple(c, A, b, K, x0, opts);
\end{lstlisting}
Note that using the simple interface, the fourth argument \mcode{x0} may be replaced by \mcode{[]}, in which case \alfonso{} will choose the default value (in this example, the all-ones vector).

\deletethis{
\subsubsection{Portfolio optimization with market impact}\label{sec:portfolio_example}

We consider a variant of the Markowitz mean-risk portfolio optimization model with two additional features: a model of transaction costs due to the trade affecting the asset prices, also known as \emph{market impact}, \citep{Almgren2010} and an explicit factor model for the return covariance matrix \citep{Perold1984}.

Let $n$ be the number of assets considered in the optimization, $f$ be the number of factors impacting the asset returns, $\vh\in\R^n$ represent the current portfolio weights, and let $\vx\in\R^n$ be the weights of the portfolio after the optimization. A factor model for the asset covariance $\vSigma$ may be given as follows:
\[ \vSigma = \vB\Omega\vB^\T + \vDelta,\]
where $\vOmega\in\R^{f\times f}$ is the factor covariance matrix, $\vB\in\R^{n\times f}$ is the factor exposure matrix, and $\vDelta\in\R^{n \times n}$ is a diagonal matrix of specific variance matrix of asset returns. 

The market impact cost is usually modeled using a power law: $\sum_{i=1}^n \delta_i |t_i|^\beta$, where $t_i = x_i-h_i$ is the trade amount and $\delta_i$ and $\beta$ are estimated parameters (positive scalars); $\beta=5/3$ is a widely accepted value \citep{AlmgrenThumHauptmannLi2005}.

If $\vmu\in\R^n$ denotes the expected returns of the $n$ assets and $\gamma$ is the highest acceptable risk level (quantified by the standard deviation of the return), then the mean-risk optimization model with market impact can be compactly written as follows:
\begin{equation}\label{eq:portfolio}
\begin{aligned}
&\underset{\vx\in\R^n}{\text{maximize}}\quad  &&\vmu^\T \vx - \sum_{i=1}^n \delta_i|x_i-h_i|^\beta \\
                &\text{subject to}  &&\vx^\T(\vB\vOmega\vB^\T + \vDelta)\vx \leq \gamma^2\\
                &                   &&\vone^\T\vx = 1\\
                &                   &&\vx \geq 0
\end{aligned}
\end{equation}
When $\beta > 1$, we can rewrite this in conic standard form by introducing the trade variables $\vt\in\R^n$ and representing the nonlinear terms in the objective using additional auxiliary variables and power cone constraints noticing that from the definition \eqref{eq:power-cone-def} we have 
\[ |t_i|^\beta \leq z_i  \iff (z_i, 1, t_i) \in \cP_{(1/\beta,\,1-1/\beta)}\]
while also rewriting the convex quadratic inequality as a second-order cone constraint:
\begin{equation}\label{eq:portfolio2}
\begin{aligned}
&\underset{\vx,\vz,\vt}{\text{minimize}}\quad  &&-\vmu^\T \vx + \vdelta^\T\vz \\
                &\text{subject to}  && \vone^\T\vx = 1\\
                &                   && \vx - \vt = \vh\\
                &                   && \vu - \vOmega^{1/2}\vB^\T\vx = \vzero\\
                &                   && \vd - \vDelta^{1/2}\vx = \vzero\\
                &                   && \vx \geq \vzero\\
                &                   && (z_i, 1, t_i) \in \cP_{(1/\beta,\,1-1/\beta)} \qquad i=1,\dots,n\\
                &                   && (\gamma, \vu, \vd) \in \cQ_{n+f+1}
\end{aligned}
\end{equation}
Strictly speaking, \eqref{eq:portfolio2} is not quite in standard form yet, as the constant $\gamma$ in the second-order constraint and $1$ in the power cone constraints appear on the left-hand side and render the cone constraints non-homogeneous. We can fix this by replacing each constant with an auxiliary decision variable that is constrained to have the correct value by a linear constraint.

Finally, the resulting optimization model can be cast and solved with \texttt{alfonso\_simple()} as a conic optimization problem where the cone is the product of a nonnegative orthant, $n$ three-dimensional power cones and an $(n+f+1)$-dimensional second-order cone. This cone can be represented using the following cone structure:
\begin{lstlisting}
K = cell(1,n+2);
K{1} = struct('type','lp','dim',n);
[K{2:n+1}] = deal(struct('type','gpow','dim',3,'lambda',[1/beta; 1-1/beta]));
K{n+2} = struct('type','socp','dim',n+f+1);
\end{lstlisting}
The code to generate the rest of the problem data (the $\vA$ matrix, the right-hand side vector $\vb$ and the objective vector $\vc$) is straightforward and we omit it here. The complete code can be found in \texttt{portfolio.m}.
}

\section{Numerical Illustration: Design of Experiments}\label{sec:design}

In this section we illustrate the potential benefit of customizable barrier computation for a semidefinite representable problem using the example of \emph{optimal design of experiments}, comparing the performance of \alfonso{} to SCS 2.1.1 and Mosek 9.2.16. For the sake of brevity, we shall forego the detailed description of the statistical problem in order to focus on the formulation of the relevant convex optimization problem, which is stated as follows \citep[Section 7.5]{BoydVandenberghe2004}.

In the optimal design problem, the input data is a (usually dense) matrix $\vV\in\R^{n \times p}$, and we seek a vector $\vx\in\R^p$ that solves the following optimization problem:
\begin{equation}\label{eq:design}
\begin{aligned}
&\underset{\vx\in\R^n}{\text{maximize}}\quad  &&\Phi(\vV\operatorname{diag}(\vx)\vV^\T)\\
                &\text{subject to}  &&\vone^\T\vx = 1\\
                &                   &&\vx \geq \vzero
\end{aligned}
\end{equation}
for some \emph{optimality criterion} $\Phi$ that maps positive definite matrices to reals. (Implicit is the constraint that the the argument of $\Phi$ is a positive definite matrix.) Most optimality criteria that are interesting from a statistical perspective are semidefinite representable in the sense of \citep{BenTalNemirovski2001}, implying that these problems are solvable using semidefinite programming. For example, the choice of $\Phi(\vM) = \lambda_\text{min}(\vM)$ leads to an \emph{E-optimal} design; see \citep[Section 7.5.2]{BoydVandenberghe2004} \deletethis{or \citep[Chapter 4]{Pukelsheim1993} }for a statistical interpretation.

Lower bound constraints on the smallest eigenvalue of a matrix $M$ can be cast in terms of a linear matrix inequality using the fact that 
$t \leq \lambda_\text{min}(\vM)$ if and only if $\vM \succcurlyeq t \vI_n$.
Since in our application $\vM = \vV\operatorname{diag}(\vx)\vV^\T$ is a linear function of our decision variables $\vx$, Eq.~\eqref{eq:design} can be readily translated to a semidefinite program with the help of an additional decision variable $t$ and subsequently solved by any semidefinite programming solver. This is the formulation that we use with Mosek and SCS. 

Instead of this semidefinite programming approach, \alfonso, equipped with a custom barrier oracle implementation, can be used to solve Eq.~\eqref{eq:design} directly as an optimization problem over the non-symmetric cone
\[ K_\vV \defeq \left\{ (t,\vx) \in \R \times \R_+^n \,\middle|\, t \leq \lambda_\text{min}(\vV\operatorname{diag}(\vx)\vV^\T)\right\}. \]
Figure \ref{fig:ex_oracle_design} shows our implementation of the $n$-LHSCB for this cone inherited from the semidefinite formulation. This example is also included with \alfonso{} in the file \texttt{examples/exp\_design/e\_design.m}; it has been slightly reformatted here to fit the page.

Table \ref{tbl:design} shows the numerical results from a set of synthetic instances of \eqref{eq:design} with $p=2n$ and $n \in \{50, 200, \dots, 500\}$, using randomly generated matrices $\vV$. Mosek and SCS were interfaced via Matlab. All computational results were obtained on a standard desktop computer equipped
with 32GB RAM and a 4 GHz Intel Core i7 processor with 4 cores running using Matlab R2017b for Windows 10. \alfonso's optimality tolerance was reduced to $\epsilon = 10^{-8}$ from the default $10^{-6}$ to match the accuracy of Mosek's solutions. Mosek and SCS were run using their default options except for increasing the maximum number of iterations for SCS to avoid early termination, tacitly acknowledging that as a first-order method, SCS is designed and expected to yield solutions with substantially lower accuracy than the interior-point methods. The solutions returned by SCS with its default tolerance settings correspond to $\epsilon\approx 10^{-3}$ in our stopping criterion. 
The complete code of this example can be found in \texttt{e\_design.m}.

In spite of returning lower-accuracy solutions, SCS exceeded one hour in the solution of the larger problems. \alfonso{} was significantly faster than both Mosek and SCS.

\afterpage{\clearpage}

\begin{table}
\caption{Solver statistics (number of iterations and total solver time in seconds) from \alfonso, Mosek 9 and SCS 2 solving the E-optimal design problem \eqref{eq:design}. \alfonso{} and Mosek returned solutions with tolerance $\epsilon \approx 10^{-8}$, the accuracy of the SCS solutions is $\epsilon \approx 10^{-3}$. Missing values indicate that the solver exceeded 1 hour.
}
	\label{tbl:design}
	\centering
	\begin{tabular}{rrrrrrr}
		\toprule
		\multicolumn{1}{c}{\multirow{2}{*}{$n$}} & \multicolumn{2}{c}{alfonso} & \multicolumn{2}{c}{Mosek} & \multicolumn{2}{c}{SCS}\\ 
\cmidrule(rl){2-3} 	\cmidrule(rl){4-5} \cmidrule(rl){6-7}
&\ \ iter\ \ & \ \ time\ \ & \ \ iter\ \ & \ \ time\ \ & \ \ iter \ \ & \ \ time\ \ \\
\midrule
50   &  48  &  0.46  & 10 &    0.52 &  3080 &    2.30  \\ 
100  &  55  &  1.37  & 11 &    2.80 &  9340 &   45.23  \\
150  &  49  &  1.51  & 11 &   12.24 & 18800 &  270.51  \\
200  &  46  &  2.15  & 11 &   33.71 & 20540 &  640.21  \\ 
250  &  51  &  5.36  & 13 &   92.88 & 40320 & 2387.81  \\              
300  &  44  &  4.41  & 10 &  155.84 & \multicolumn{2}{c}{$> 1$ hr}\\
350  &  50  &  9.25  & 11 &  304.98 &   \\
400  &  46  &  7.94  & 11 &  521.45 &   \\
450  &  57  & 20.70  & 12 &  908.56 &   \\
500  &  51  & 12.85  & 12 & 1420.50 &   \\
	\bottomrule\\
	\end{tabular}
\end{table}

\afterpage{\clearpage}

\begin{figure}
\caption{A membership and barrier function oracle for the E-optimal design example.}\label{fig:ex_oracle_design}
\begin{lstlisting}[aboveskip=1ex,belowskip=0ex]
function [in, g, H, L] = e_design(tx, pars)
% This function implements a membership and barrier function oracle for
% the E-optimal design example e_design.m
%
% INPUT
% tx:                   column vector representing [t; x(1); ...; x(n)]
% pars:                 structure with a single field pars.v
%                       pars.v is a two-dimensional array whose ith column
%                       v(:,i) represents the ith design vector (i=1,...,p)
t = tx(1);
x = tx(2:end);
[n,p] = size(pars.v);
 
% in the cone?
if any(x <= 0)
    in = false;  g = NaN;  H = NaN;  L = NaN;  return
end
 
Ax = -t*eye(n) + pars.v*diag(x)*pars.v';
 
[L,err] = chol(Ax,'lower');
if err > 0
    in = false;  g = NaN;  H = NaN;  L = NaN;  return
else
    in = true;
end
 
% compute g and H if required
if nargout > 1
    g = [0; -1./x];
    
    Li = inv(L);
    g(1) = Li(:)'*Li(:);
    w = L\pars.v;
    for i=1:p
       g(i+1) = -w(:,i)'*w(:,i);
    end
    
    % compute H and L if required
    if nargout > 2
        H = diag([0; x.^(-2)]);
        
        invAx = Li'*Li;
        H(1,1) = H(1,1) + invAx(:)'*invAx(:);
        Lws = L' \ w;
        for i=2:p+1
            H(i,1) = H(i,1) - Lws(:,i-1)'*Lws(:,i-1);
        end
        H(1,2:p+1) = H(2:p+1,1)';
        H(2:end,2:end) = H(2:end,2:end) + (w'*w).^2;
        
        if nargout > 3
            [L,err] = chol(H,'lower');
            if err > 0
                in = false;  g = NaN;  H = NaN;  L = NaN;  return
            end
        end
    end
end
return
\end{lstlisting}
\end{figure}

\afterpage{\clearpage}

\section{Discussion}

\alfonso{} provides an easily usable and customizable, yet efficient, open-source tool for conic optimization. Using its oracle interface, researchers and practitioners can solve optimization problems over non-symmetric cones that do not have a convenient representation in terms of symmetric cone constraints. Additionally, as our last example shows, it can even provide a significant speedup over state-of-the-art solvers in problems with a straightforward semidefinite programming formulation by exploiting problem structure and avoiding the introduction of a large number of auxiliary variables. A key feature of the underlying algorithm is that all of its parameters are generic, applicable to any convex cone. Therefore the user only needs to provide the code to compute the derivatives of the barrier function and a point in the interior of the cone.

\paragraph{Extending the simple interface} The simple interface currently supports a limited number of non-symmetric cones (mostly the same ones as SCS and Mosek). New cones can be easily added with minimal changes to the code, limited to a single file \texttt{src/alfonso\_simple.m}. Specifically, once the membership and barrier function oracle is prepared (as a separate Matlab file), the simple interface only needs a pointer to the cone and an interior point, both added in the form of a new line in a \mcode{switch-case} structure.

\deletethis{
\paragraph{Limitations} As with all primal-dual interior-point methods using a homogeneous self-dual embedding, the algorithm is guaranteed to find a primal-dual optimal solution pair or a certificate of infeasibility only when the problem satisfies strong duality. There has been some recent work indicating that when self-dual IPMs fail to certify optimality or infeasibility (because in the limit both $\tau$ and $\kappa$ tend to $0$), the output may be used as a facial reduction certificate \citep{PermenterFribergAndersen2017}. In its current version, \alfonso{} does not implement any facial reduction.

\paragraph{Efficiency} The theoretical worst-case iteration complexity of the algorithm, $\Oh(\sqrt{\nu}\log(1/\epsilon))$, is of the same order of magnitude as the iteration complexity of popular symmetric cone optimization methods. The bulk of the work in each iteration is the solution of the Newton systems in the predictor and corrector steps. This defines a practical bound on what an LHSCB with ``efficiently computable'' gradient and Hessian means. Ideally, the computation of the derivatives should be faster than the solution of the Newton steps, which takes $\Oh(n^3)$ flops; this ensures that the barrier function oracle does not become the bottleneck.
}

%
%
%


\newpage
\bibliographystyle{informs2014} 
\bibliography{alfonso_ijoc} 

\begin{thebibliography}{21}
\providecommand{\natexlab}[1]{#1}
\providecommand{\url}[1]{\texttt{#1}}
\providecommand{\urlprefix}{URL }

\bibitem[{Andersen(2002)}]{Andersen2002}
Andersen ED (2002) Handling free variables in primal-dual interior-point
  methods using a quadratic cone. \emph{Proceedings of the SIAM Conference on
  Optimization, Toronto}.

\bibitem[{Ben-Tal \protect\BIBand{} Nemirovski(2001)}]{BenTalNemirovski2001}
Ben-Tal A, Nemirovski A (2001) \emph{Lectures on Modern Convex Optimization}
  (Philadelphia, {PA}: {SIAM}), ISBN 0-89871-491-5.

\bibitem[{Blekherman et~al.(2013)Blekherman, Parrilo, \protect\BIBand{}
  Thomas}]{BlekhermanParriloThomas2013}
Blekherman G, Parrilo PA, Thomas RR, eds. (2013) \emph{Semidefinite
  optimization and convex algebraic geometry} (Philadelphia, PA: {SIAM}), ISBN
  978-1-611972-28-3.

\bibitem[{Boyd \protect\BIBand{} Vandenberghe(2004)}]{BoydVandenberghe2004}
Boyd SP, Vandenberghe L (2004) \emph{Convex Optimization} (Cambridge University
  Press), ISBN 0-521-83378-7.

\bibitem[{Chares(2009)}]{Chares2009}
Chares R (2009) \emph{Cones and interior-point algorithms for structured convex
  optimization involving powers and exponentials}. Ph.D. thesis, Universit{\'e}
  Catholique de Louvain.

\bibitem[{Coey et~al.(2020)Coey, Kapelevich, \protect\BIBand{}
  Vielma}]{CoeyKapelevichVielma2020}
Coey C, Kapelevich L, Vielma JP (2020) Towards practical generic conic
  optimization. \emph{arXiv preprint arXiv:2005.01136} .

\bibitem[{Domahidi et~al.(2013)Domahidi, Chu, \protect\BIBand{} Boyd}]{ecos}
Domahidi A, Chu E, Boyd S (2013) {ECOS}: {A}n {SOCP} solver for embedded
  systems. \emph{European Control Conference (ECC)}, 3071--3076.

\bibitem[{Glineur \protect\BIBand{} Terlaky(2004)}]{GlineurTerlaky2004}
Glineur F, Terlaky T (2004) Conic formulation for l$_p$-norm optimization.
  \emph{Journal of Optimization Theory and Applications} 122(2):285--307.

\bibitem[{G{\"u}ler(1997)}]{Guler1997}
G{\"u}ler O (1997) Hyperbolic polynomials and interior point methods for convex
  programming. \emph{Mathematics of Operations Research} 22(2):350--377.

\bibitem[{Iliman \protect\BIBand{} de~Wolff(2016)}]{IlimanDeWolff2016}
Iliman S, de~Wolff T (2016) Amoebas, nonnegative polynomials and sums of
  squares supported on circuits. \emph{Research in the Mathematical Sciences}
  3(1):9, \urlprefix\url{http://dx.doi.org/10.1186/s40687-016-0052-2}.

\bibitem[{Karimi \protect\BIBand{} Tun\c{c}el(2019)}]{KarimiTuncel2019}
Karimi M, Tun\c{c}el L (2019) Domain-driven solver ({DDS}): a {MATLAB}-based
  software package for convex optimization problems in domain-driven form.
  \emph{arXiv preprint arXiv:1908.03075} .

\bibitem[{{MOSEK ApS}(2019)}]{mosek9}
{MOSEK ApS} (2019) {MOSEK} {O}ptimization {S}uite release {9.0.105}.
  \urlprefix\url{https://docs.mosek.com/9.0/releasenotes/index.html}.

\bibitem[{Nesterov \protect\BIBand{}
  Nemirovskii(1994)}]{NesterovNemirovskii1994}
Nesterov Y, Nemirovskii A (1994) \emph{Interior-point polynomial algorithms in
  convex programming}, volume~13 of \emph{SIAM Studies in Applied Mathematics}
  (Philadelphia, PA: Society for Industrial and Applied Mathematics (SIAM)),
  ISBN 0-89871-319-6,
  \urlprefix\url{http://dx.doi.org/10.1137/1.9781611970791}.

\bibitem[{O'Donoghue et~al.(2016)O'Donoghue, Chu, Parikh, \protect\BIBand{}
  Boyd}]{scs}
O'Donoghue B, Chu E, Parikh N, Boyd S (2016) Operator splitting for conic
  optimization via homogeneous self-dual embedding. \emph{Journal of
  Optimization Theory and Applications} 169:1042–1068,
  \urlprefix\url{http://dx.doi.org/10.1007/s10957-016-0892-3}.

\bibitem[{Papp(2019)}]{Papp2019}
Papp D (2019) Duality of sum of nonnegative circuit polynomials and optimal
  {SONC} bounds. \emph{arXiv preprint arXiv:1912.04718}
  \urlprefix\url{https://arxiv.org/abs/1912.04718}.

\bibitem[{Papp \protect\BIBand{}
  Y{\i}ld{\i}z(2017)}]{PappYildiz2017corrigendum}
Papp D, Y{\i}ld{\i}z S (2017) On ``{A} homogeneous interior-point algorithm for
  non-symmetric convex conic optimization". \emph{arXiv preprint
  arXiv:1712.00492} \urlprefix\url{https://arxiv.org/abs/1712.00492}.

\bibitem[{Papp \protect\BIBand{} Y{\i}ld{\i}z(2019)}]{PappYildiz2019}
Papp D, Y{\i}ld{\i}z S (2019) Sum-of-squares optimization without semidefinite
  programming. \emph{SIAM Journal on Optimization} 29(1):822--851,
  \urlprefix\url{http://dx.doi.org/10.1137/17M1160124}.

\bibitem[{Renegar(2001)}]{Renegar2001}
Renegar J (2001) \emph{A mathematical view of interior-point methods in convex
  optimization}. MOS-SIAM Series on Optimization (Phiadelphia, PA: Society for
  Industrial and Applied Mathematics (SIAM)), ISBN 0-89871-502-4,
  \urlprefix\url{http://dx.doi.org/10.1137/1.9780898718812}.

\bibitem[{Renegar(2004)}]{Renegar2004}
Renegar J (2004) Hyperbolic programs, and their derivative relaxations.
  Technical report, Cornell University Operations Research and Industrial
  Engineering.

\bibitem[{Roy \protect\BIBand{} Xiao(2018)}]{RoyXiao2018}
Roy S, Xiao L (2018) On self-concordant barriers for generalized power cones.
  Technical Report MSR-TR-2018-3, Microsoft Research,
  \urlprefix\url{https://www.microsoft.com/en-us/research/publication/on-self-concordant-barriers-for-generalized-power-cones/}.

\bibitem[{Skajaa \protect\BIBand{} Ye(2015)}]{SkajaaYe2015}
Skajaa A, Ye Y (2015) A homogeneous interior-point algorithm for nonsymmetric
  convex conic optimization. \emph{Mathematical Programming Ser. A}
  150(2):391--422, ISSN 0025-5610,
  \urlprefix\url{http://dx.doi.org/10.1007/s10107-014-0773-1}.

\end{thebibliography}


\end{document}